\documentclass[a4paper]{amsart}

\usepackage[utf8]{inputenc}
\usepackage[english]{babel}

\usepackage{enumerate}
\usepackage[usenames,dvipsnames,table]{xcolor}
\usepackage{color}
\usepackage{amsmath,amsthm,amssymb}
\usepackage{amsfonts}

\usepackage{fourier}

\usepackage{todonotes}
\usepackage{mathrsfs}

\usepackage{hyperref}
\definecolor{colorcita}{RGB}{21,86,130}
\definecolor{colorref}{RGB}{5,10,177}
\definecolor{colorweb}{RGB}{177,6,38}

\newcommand{\N}{\mathbb{N}}
\newcommand{\C}{\mathbb{C}}
\newcommand{\R}{\mathbb{R}}

\newcommand{\E}{\mathscr{E}} 
\newcommand{\D}{\mathscr{D}}
\newcommand{\supp}{\mbox{supp}\,}

\newtheorem{theorem}{Theorem}[section]
\newtheorem{lemma}[theorem]{Lemma}
\newtheorem{corollary}[theorem]{Corollary}
\newtheorem{proposition}[theorem]{Proposition}

\theoremstyle{definition}
\newtheorem{remark}[theorem]{Remark}

\newtheorem{example}[theorem]{Example}

\newcommand{\sign}{\mbox{sign}}

\title{Mean ergodic composition operators on spaces of smooth functions and distributions}

\author[T.\ Kalmes]{Thomas Kalmes}
\address{Thomas Kalmes, Chemnitz University of Technology, Faculty of Mathematics, 09107 Chemnitz, Germany}
\email{thomas.kalmes@math.tu-chemnitz.de}

\author[D.\ Santacreu]{Daniel Santacreu}
\address{Daniel Santacreu, Instituto Universitario Matem\'atica Pura y Aplicada IUMPA, Universitat Polit\`ecnica de Val\`encia, Camino de Vera, s/n, 46701 Valencia, Spain}
\email{dasanfe5@posgrado.upv.es}

\thanks{This manuscript version is made available under the CC-BY-NC-ND 4.0 license \href{https://creativecommons.org/licenses/by-nc-nd/4.0/}{http://creativecommons.org/licenses/by-nc-nd/4.0/}\newline First published in Proceedings of the American Mathematical Society in vol.\ 150(6), 2022, pages 2603-2616, published by American Mathematical Society\newline
	DOI: 10.1090/proc/15894\newline
	\href{https://doi.org/10.1090/proc/15894}{https://doi.org/10.1090/proc/15894}}

\begin{document}
\maketitle

\begin{abstract}
	We investigate (uniform) mean ergodicity of weighted composition operators on the space of smooth functions and the space of distributions, both over an open subset of the real line. Among other things, we prove that a composition operator with a real analytic diffeomorphic symbol is mean ergodic on the space of distributions if and only if it is periodic with period 2. Our results are based on a characterization of mean ergodicity in terms of Ces\`aro boundedness and a growth property of the orbits for operators on Montel spaces which is of independent interest.\\
	
	
\end{abstract}

\section{Introduction}

In this note we contribute to the investigation of the dynamical behaviour of weighted composition operators $C_{w,\phi}(u)=w\cdot (u\circ\phi)$ on the space of smooth (i.e.\ infinitely many times differentiable) functions $\E(X)$ on an open subset $X\subseteq\R$ as well as on the space of distributions $\mathscr{D}'(X)$, where the symbol $\phi$ of $C_{w,\phi}$ is a smooth self map of $X$, respectively a diffeomorphism of $X$ when dealing with $\mathscr{D}'(X)$, and the weight $w$ is a complex valued smooth function on $X$. The space of smooth functions $\E(X)$ is endowed with its natural topology of uniform convergence on compact subsets of $X$ of all derivatives up to an arbitrary finite order while $\mathscr{D}'(X)$ is equipped with its strong dual topology, being the topological dual space of $\mathscr{D}(X)$, the space of test functions on $X$. Thus, both spaces are Montel spaces. We are interested in when such weighted composition operators are mean ergodic (definitions will be given in Section \ref{sec: general abstract results} below).

In recent years there have been several articles studying mean ergodicity and related properties of (weighted) composition operators on various spaces of functions, such as spaces of holomorphic functions in finite dimensions \cite{BoDo2011A}, \cite{BGJJ2016mw}, \cite{BeGCJoJo16},\cite{HaShZh19}, \cite{Beltran20}, \cite{JoRA20}, \cite{SeMeBo20}, spaces of holomorphic functions on infinite dimensional Banach spaces \cite{JoSaSP21}, spaces of homogeneous polynomials on infinite dimensional Banach spaces \cite{JoSaSP20}, spaces of real analytic functions \cite{BoDo11B}, the Schwartz space of rapidly decreasing functions on $\R$ \cite{FeGaJo18}, spaces of meromorphic functions \cite{GCJoJo16}, and within the general framework of function spaces defined by local properties \cite{K2019p}.  

This note is organized as follows. In Section \ref{sec: general abstract results} we show that for a continuous linear operator on a Montel space the properties of mean ergodicity and uniform mean ergodicity coincide, and we give a characterization of these properties in terms of Ces\`aro boundedness of the operator together with a growth property of its orbits (Theorem \ref{theo:ergodicity and transposed} (b)). In Section \ref{sec:weighted on smooth}, based on the aforementioned result, we derive necessary and sufficient conditions for $C_{w,\phi}$ to be mean ergodic on $\E(X)$ (Theorem \ref{theo: mean ergodicity on smooth}). Under the additional assumption that $\phi$ is a diffeomorphism and $\{x\in X;\,w(x)\neq 0\}$ is dense in $X$, in Section \ref{sec: mean ergodic on D'} we show that mean ergodicity of $C_{w,\phi}$ on $\mathscr{D}'(X)$ forces a rather restrictive behaviour of the symbol $\phi$, namely $\phi$ as well as $\phi^{-1}$ have stable orbits (Theorem \ref{theo: both stable orbits}). This restrictive property is then used to show that for a real analytic diffeomorphism $\phi$ the corresponding unweighted composition operator $C_\phi:=C_{1,\phi}$ is mean ergodic on $\mathscr{D}'(X)$ if and only $C_\phi$ is periodic of period 2 (Theorem \ref{theo: real analytic diffeos}).

\section{General abstract results}\label{sec: general abstract results}

Let $E$ be a locally convex Hausdorff space (briefly, lcHs) and $T\in \mathcal{L}(E)$, where as usual we denote by $\mathcal{L}(E)$ the space of continuous linear operators on $E$. Moreover, by $cs(E)$ we denote the set of continuous seminorms on $E$. $T$ is said to be {\it topologizable} if for every $p\in cs(E)$ there is $q\in cs(E)$ such that for every $m\in\N$ there is $\gamma_m>0$ with
\[p\left(T^mx\right)\leq \gamma_m q(x) \text{ for all } x\in E.\]
For the special case that in the above inequality one can take $\gamma_m=1$ for all $m\in\N$ we say that $T$ is {\it power bounded}. In this case the family $\{T^m:\,m\in\N\}$ is an equicontinuous subset of $\mathcal{L}(E)$. Moreover, $T$ is {\it Cesàro bounded} if the family $\{T^{[n]}:\,n\in\N\}$ is an equicontinuous subset of $\mathcal{L}(E)$, where $T^{[n]}$ denotes the $n$-th Cesàro mean given by
\begin{equation*}
	\frac{1}{n}\sum_{m=1}^{n} T^m
\end{equation*}
An operator $T\in\mathcal{L}(E)$ is called {\it mean ergodic} if there is $P\in\mathcal{L}(E)$ such that  for each $x\in E$ it holds $\lim_{n\to\infty}T^{[n]}x=Px$. In case that the convergence is uniform on bounded subsets of $E$ then $T$ is called {\it uniformly mean ergodic}. For $T\in \mathcal{L}(E)$ and $n\in\N$ we have the following identities (where $T^{[0]}=I$)
\begin{equation}\label{ME implies second}
	\frac{1}{n}T^{n}=T^{[n]}-\frac{n-1}{n}T^{[n-1]},
\end{equation}
\begin{equation}\label{I-T means}
	(I-T)T^{[n]}=T^{[n]}(I-T)=\frac{1}{n}(T-T^{n+1})
\end{equation}
so that $\lim_{n\to\infty} \frac{1}{n}T^nx=0$ for every $x\in E$ whenever $T$ is mean ergodic. 

The following theorem is a special case of Eberlein's mean ergodic theorem which is proved by a straightforward modification of the proof in \cite[Chapter 2, § 2.1, Theorem 1.5, p.\ 76]{Krengel}. In our context, one has to set the semigroup of operators $\mathcal{S}=\{T^n:\,n\in\N_0\}$ and the ergodic net $\{T^{[n]}:\,n\in\N\}$, where one has to take into account the fact that due to
$$T^k T^{[n]} x- T^{[n]} x=T^{[n]} T^k x- T^{[n]} x = \frac{1}{n} \sum_{m=1}^{k} T^{n+m} x - \frac{1}{n} \sum_{m=1}^{k} T^{m} x$$
it follows $\lim_{n\rightarrow\infty}T^{[n]} T^k x- T^{[n]} x=\lim_{n\rightarrow\infty}T^k T^{[n]} x- T^{[n]} x=0$ in $E$ whenever $x\in E$ satisfies
\begin{equation}\label{ME second property}
	\lim_{n\to\infty} \frac{1}{n}T^nx=0.
\end{equation}

\begin{theorem}\label{Eberlein}
	Let $E$ be a lcHS and let $T\in \mathcal{L}(E)$ be Cesàro bounded and let $x\in E$ be such that $\lim_{n\to\infty} \frac{1}{n}T^nx=0$. The following conditions are equivalent for $y\in E$:
	\begin{enumerate}
		\item[(a)] $Ty=y$ and $y$ belongs to the closed convex hull of the orbit $O(x,T):=\{T^m x:\, m\in\N_0\}$ of $x$. 
		\item[(b)] $y=\lim\limits_{n\to\infty}T^{[n]}x$
		\item[(c)] $y=\sigma(E,E^\prime)-\lim\limits_{n\to\infty}T^{[n]}x$
		\item[(d)] $y$ is a $\sigma(E,E^\prime)$-cluster point of $\left(T^{[n]}x\right)_{n\in\N}$, i.e.\ for each $0$-neighborhood $U$ in $(E,\sigma(E,E^\prime))$ and each $m\in\N$ there is $n>m$ such that $y-T^{[n]}x\in U$.
	\end{enumerate}
\end{theorem}

As a consequence of Theorem \ref{Eberlein} we obtain the next result.

\begin{corollary}\label{ME coro}
	Let $T\in \mathcal{L}(E)$ be Cesàro bounded. Then $T$ is mean ergodic if and only if \eqref{ME second property} is satisfied for all $x\in E$ and $\left(T^{[n]}x\right)_{n\in\N}$ is relatively $\sigma(E,E^\prime)$-compact for each $x\in E$.
\end{corollary}
\begin{proof}
	If $T$ is mean ergodic, \eqref{ME second property} holds by identity \eqref{ME implies second}. Clearly $\left(T^{[n]}x\right)_{n\in\N}$ is relatively $\sigma(E,E^\prime)$-compact for each $x\in E$ because $\left(T^{[n]}x\right)_{n\in\N}$ converges in $E$. 
	
	Conversely, given $x\in E$, the set $\left\{T^{[n]}x:\,n\in\N\right\}$ is relatively $\sigma(E,E^\prime)$-compact, and therefore $\left(T^{[n]}x\right)_{n\in\N}$ has a  $\sigma(E,E^\prime)$-cluster point $y\in E$. By Theorem \ref{Eberlein}, necessarily $y=\lim\limits_{n\to\infty}T^{[n]}x$. We define 
	\[Px:=\lim\limits_{n\to\infty}T^{[n]}x\]
	for each $x\in E$. Since $\left(T^{[n]}\right)_{n\in\N}$ is equicontinuous in $\mathcal{L}(E)$ we obtain that $P\in\mathcal{L}(E)$ and $T$ is mean ergodic.
\end{proof}

The following result is a version of \cite[Theorem 2.4]{ABR2009} for reflexive locally convex spaces.

\begin{theorem}\label{mean ergodic}
	Let $E$ be a reflexive lcHs and $T\in \mathcal{L}(E)$. Then $T$ is mean ergodic if and only if $T$ is Cesàro bounded and $\lim_{n\to\infty} \frac{1}{n}T^n x=0$ for every $x\in E$.
\end{theorem}
\begin{proof}
	We recall that a locally convex Hausdorff space $E$ is reflexive if and only if it is barrelled and every bounded subset of $E$ is relatively $\sigma(E,E^\prime)$-compact.
	
	As in the proof of Corollary \ref{ME coro}, for mean ergodic $T$ we have that \eqref{ME second property} is satisfied for all $x\in E$ and $\{T^{[n]}x:\,n\in\N\}$ is relatively $\sigma(E,E^\prime)$-compact for each $x\in E$. Consequently, $\{T^{[n]};n\in\N\}\subset \mathcal{L}(E)$ is a pointwise bounded set and by the Uniform Boundedness Principle $T$ is Cesàro bounded.
	
	Conversely, if $\{T^{[n]}:\,n\in\N\}$ is equicontinuous, for each $x\in E$ the set $\{T^{[n]}x;n\in\N\}$ is bounded. Thus, as $E$ is reflexive, this set is relatively $\sigma(E,E^\prime)$-compact. An application of Corollary \ref{ME coro} concludes the proof.
\end{proof}

Since in semi-reflexive spaces bounded subsets are relatively $\sigma(E,E^\prime)$-compact and since for power bounded $T$ condition \eqref{ME second property} trivially holds, Corollary \ref{ME coro} yields the next result from \cite[Proposition 3.3]{BoPaRi11}.

\begin{corollary}\label{coro: power bounded on semi-reflexive are mean ergodic}
	Every power bounded operator on a semi-reflexive locally convex Hausdorff space is mean ergodic. 
\end{corollary}

We are now ready to prove the following result which contains a characterization of (uniform) mean ergodicity for an operator on Montel spaces in terms of Ces\`aro boundedness and a growth property of its orbits. 

\begin{theorem}\label{theo:ergodicity and transposed}
	Let $E$ be a Montel space and let $T\in \mathcal{L}(E)$.
	\begin{enumerate}
		\item[(a)] $T$ is mean ergodic if and only if $T$ is uniformly mean ergodic.
		\item[(b)] The following are equivalent.
		\begin{enumerate}
			\item[(i)] $T$ is Ces\`aro bounded and $\lim_{n\rightarrow\infty}\frac{T^n}{n}=0$, pointwise in $E$.
			\item[(ii)] $T$ is mean ergodic on $E$.
			\item[(iii)] $T$ is uniformly mean ergodic on $E$.
			\item[(iv)] $T^t$ is mean ergodic on $(E',\beta(E',E))$.
			\item[(v)] $T^t$ is uniformly mean ergodic on $(E',\beta(E',E))$.
			\item[(vi)] $T^t$ is Ces\`aro bounded on $(E',\beta(E',E))$ and $\lim_{n\rightarrow\infty}\frac{(T^{t})^n}{n}=0$, pointwise in $(E',\beta(E',E))$.
		\end{enumerate}
	\end{enumerate}
\end{theorem}

\begin{proof}
	
	Trivially, every uniformly mean ergodic operator is mean ergodic. Let $T$ be mean ergodic. Then, $\left\{T^{[n]}:\,n\in\N\right\}$ is equicontinuous and since $E$ is a Montel space, every bounded subset $B$ of $E$ is relatively compact. Since on equicontinuous subsets of $\mathcal{L}(E)$ pointwise convergence on $E$ and uniform convergence of relatively compact subsets of $E$ coincide, it follows that $\left(T^{[n]}\right)_{n\in\N}$ converges uniformly on bounded subsets of $E$. Thus, (a) is proved.
	
	In order to prove (b), we observe that (i) and (ii) as well as (iv) and (vi) are equivalent by Theorem \ref{mean ergodic} while (ii) and (iii) are equivalent by part (a). Since with $E$ also $(E',\beta(E',E))$ is a Montel space, the equivalence of (iv) and (v) follows from part (a) as well. Finally, by \cite[Corollary 2.7 (ii)]{ABR2010} and the fact that Montel spaces are reflexive, (iii) and (v) are equivalent. 
\end{proof}

\section{Weighted composition operators on $\E(X)$}\label{sec:weighted on smooth}

In this section we study the mean ergodicity of weighted composition operators $C_{w,\phi}$ on the space of smooth functions $\E(X)$, where $X\subseteq \R$ is an open set. Here, $w:X\rightarrow\C$ and $\phi:X\rightarrow X$ are smooth functions and $\E(X)$ is equipped with its standard topology, i.e.\ with the Fr\'echet space topology generated by the seminorms
\begin{equation*}
	\forall\,K\subset X\mbox{ compact}, s\in\N_0, f\in \E(X):\,\|f\|_{s,K}:=\sup_{x\in K, 0\leq r\leq s}\left|f^{(r)}(x)\right|.
\end{equation*}
As usual, $C_{w,\phi}:\E(X)\rightarrow\E(X)$ is defined as
\begin{equation*}
	C_{w,\phi} f(x):= w(x) f(\phi(x)),\, x\in X,
\end{equation*}
for all $f\in\E(X)$. Then, $\E(X)$ is a nuclear Fr\'echet space and thus, in particular a Montel space and clearly $C_{w,\phi}\in\mathcal{L}\left(\mathscr{E}(X)\right)$. Thus, by Theorem \ref{theo:ergodicity and transposed} (b), $C_{w,\phi}$ is mean ergodic if and only if it is uniformly mean ergodic.

For $s\in\N$ we have, using Leibniz' rule and Faà di Bruno's formula \cite{Johnson02},
\begin{align*}
	\left(C_{w,\phi}^n f\right)^{(s)}=&\left(\prod_{l=0}^{n-1}w(\phi^l)\right)^{(s)} f(\phi^n)\\
	&+\sum_{r=1}^s \binom{s}{r}\left(\prod_{l=0}^{n-1}w(\phi^l)\right)^{(s-r)}\sum_{j=1}^r f^{(j)}(\phi^n)B_{r,j}((\phi^n)',\ldots,(\phi^n)^{(r-j+1)}),
\end{align*}
where $B_{r,j}$ denote the corresponding Bell polynomials. Please note that by $\phi^l$ we denote the $l$-fold composition of $\phi$ with itself etc.. In order to simplify our notation, we abbreviate
\begin{equation*}
	\forall\,r\in\N, j\in\{1,\ldots,r\}, n\in\N_0: B^\phi_{r,j,n}:=B_{r,j}((\phi^n)',\ldots,(\phi^n)^{(r-j+1)})
\end{equation*}
as well as
\begin{equation*}
	\forall\,r\in\N, n\in\N_0, x\in \R: B^\phi_{0,0,n}(x):=1, B^\phi_{r,0,n}(x):=0.
\end{equation*}
With this notation, we have
\begin{equation}\label{eq:s-th derivative}
	\left(C_{w,\phi}^n f\right)^{(s)}=\sum_{0\leq j\leq r\leq s}\binom{s}{r}\left(\prod_{l=0}^{n-1} w(\phi^l)\right)^{(s-r)} f^{(j)}(\phi^n) B^\phi_{r,j,n}
\end{equation}
for every $f\in\E(X)$ and $s\in \N_0$. Evaluating this equality for the special case of $f_\lambda(y):=\exp(\lambda y), \lambda\in \C, y\in\R,$ yields
\begin{equation}\label{eq:s-th derivative of exponential}
	\forall\,\lambda\in\C, s\in\N_0:\,\left(C_{w,\phi}^n f_\lambda\right)^{(s)}=f_\lambda\circ\phi^n\sum_{0\leq j\leq r\leq s}\binom{s}{r}\left(\prod_{l=0}^{n-1} w(\phi^l)\right)^{(s-r)} \lambda^j B^\phi_{r,j,n}.
\end{equation}

Now we discuss necessary and sufficient conditions involving mean ergodicity. The following result is a characterization of the property \eqref{ME second property} for $C_{w,\phi}$. Recall, that $\phi:X\rightarrow X$ is said to have \emph{stable orbits} if for each compact $K\subset X$ there is another compact subset $L\subset X$ with $\phi^n(K)\subseteq L$ for every $n\in\N_0$.
\begin{proposition}\label{prop:Bell to 0}
	Let $\phi:X\rightarrow X$ and $w:X\rightarrow\C$ be smooth functions such that $\{x\in X;\,w(\phi^m(x))\neq 0\}$ is dense in $X$ for every $m\in\N_0$. Then, the following are equivalent:
	\begin{enumerate}
		\item[(i)] $\lim_{n\rightarrow\infty}\frac{1}{n}C_{w,\phi}^n f=0$ in $\E(X)$.
		\item[(ii)] $\phi$ has stable orbits and for every compact $K\subset X$, $s\in\N_0$, and $h\in\{0,\ldots,s\}$ it holds
		\begin{equation*}
			\lim_{n\rightarrow\infty}\frac{1}{n}\left\|\sum_{r=h}^s\left(\prod_{l=0}^{n-1} w(\phi^l)\right)^{(s-r)}B^\phi_{r,h,n}\right\|_{0,K}=0.
		\end{equation*}
	\end{enumerate}
\end{proposition}
\begin{proof}
	In order to show that (i) implies (ii), for $s\in\N$ we set $\lambda_s:=\exp(i\frac{2\pi}{s})$ and for $h\in\{1,\ldots,s\}$ we define
	\begin{equation*}
		Q_{s,h}(x):=\prod_{1\leq j\leq s, \hspace{0.05cm} j\neq h} (\lambda_s^j-x), x\in\R.
	\end{equation*}
	Thus, $Q_{s,h}$ is a polynomial of degree $s-1$ with $Q_{s,h}(\lambda_s^h)\neq 0$. Then,
	\begin{equation*}
		P_{s,h}(x):=\frac{1}{Q_{s,h}(\lambda_s^h)} Q_{s,h}(x), x\in\R,
	\end{equation*}
	is a polynomial of degree $s-1$ satisfying $P_{s,h}(\lambda_s^j)=\delta_{j,h}$ for $j\in\{1,\ldots,s\}$ where $\delta_{j,h}$ denotes Kronecker's delta. Let $\alpha_0^{(s,h)},\ldots, \alpha_{s-1}^{(s,h)}\in\C$ be such that $P_{s,h}(x)=\sum_{k=0}^{s-1}\alpha_k^{(s,h)}x^k$.
	
	Since $\E(X)$ is a Fr\'echet space, by the Uniform Boundedness Principle, (i) implies the equicontinuity of $\left\{\frac{1}{n}C_{w,\phi}^n;\,n\in\N_0\right\}$. In particular, $C_{w,\phi}$ is topologizable, so that by \cite[Corollary 3.12 and proof of Corollary 3.15]{K2019p} $\phi$ has stable orbits. Now we fix $s\in\N$ and a compact $K\subset X$. For arbitrary $k\in\N_0$, applying (i) and (\ref{eq:s-th derivative of exponential}) to $f_{\lambda_s^k}$ and $\lambda_s^k$, respectively, yields
	\begin{equation}\label{eq:auxiliary 1}
		0=\lim_{n\rightarrow\infty}\frac{1}{n}\left\|f_{\lambda_s^k}\left(\phi^n\right)\left(\sum_{0\leq j\leq r\leq s}\binom{s}{r}\left(\prod_{l=0}^{n-1}w\left(\phi^l\right)\right)^{(s-r)}\lambda_s^{jk}B^\phi_{r,j,n}\right)\right\|_{0,K}.
	\end{equation}
	Because $\phi$ has stable orbits, there is a compact set $L\subset X$ such that $\phi^n(K)\subseteq L$ for every $n\in\N_0$. By compactness of $L\times\{\lambda\in\C;\,|\lambda|=1\}$, there is $C>0$ such that
	\begin{equation*}
		\forall\,n,k\in\N_0, x\in K:|f_{\lambda_s^k}(\phi^n(x))|\geq C.
	\end{equation*} 
	Thus, (\ref{eq:auxiliary 1}) implies for every $h\in\{1,\ldots,s\}$
	\begin{align*}
		0&=\lim_{n\rightarrow\infty}\frac{1}{n}\left\|\sum_{k=0}^{s-1} \alpha_k^{(s,h)}\sum_{0\leq j\leq r\leq s}\binom{s}{r}\left(\prod_{l=0}^{n-1} w\left(\phi^l\right)\right)^{(s-r)}\lambda_s^{jk}B^\phi_{r,j,n}\right\|_{0,K}\\
		&=\lim_{n\rightarrow\infty}\frac{1}{n}\left\|\sum_{0\leq j\leq r\leq s}\binom{s}{r}\left(\prod_{l=0}^{n-1} w\left(\phi^l\right)\right)^{(s-r)}\left(\sum_{k=0}^{s-1} \alpha_k^{(s,h)}\lambda_s^{jk}\right)B^\phi_{r,j,n}\right\|_{0,K}\\
		&=\lim_{n\rightarrow\infty}\frac{1}{n}\left\|\sum_{0\leq j\leq r\leq s}\binom{s}{r}\left(\prod_{l=0}^{n-1} w\left(\phi^l\right)\right)^{(s-r)} P_{s,h}\left(\lambda_s^j\right)B^\phi_{r,j,n}\right\|_{0,K}\\
		&=\lim_{n\rightarrow\infty}\frac{1}{n}\left\|\sum_{r=h}^s \binom{s}{r}\left(\prod_{l=0}^{n-1} w\left(\phi^l\right)\right)^{(s-r)}B^\phi_{r,h,n}\right\|_{0,K}.
	\end{align*}
	Moreover, evaluating (i) for $f=1$ implies
	\begin{equation*}
		\lim_{n\rightarrow\infty}\left\|\left(\prod_{l=0}^{n-1}w\left(\phi^j\right)\right)^{(s)}\right\|_{0,K}=0.
	\end{equation*}
	Thus, combining the last two equalities gives
	\begin{equation}\label{eq:result for natural s}
		\forall\, h\in\{0,\ldots, s\}: 0=\lim_{n\rightarrow\infty}\frac{1}{n}\left\|\sum_{r=h}^s \binom{s}{r}\left(\prod_{l=0}^{n-1} w\left(\phi^l\right)\right)^{(s-r)}B^\phi_{r,h,n}\right\|_{0,K}.
	\end{equation}
	Finally, evaluating (i) again for $f=1$ implies for $s=0$ that
	\begin{equation*}
		0=\lim_{n\rightarrow\infty}\frac{1}{n}\left\|\sum_{r=0}^0 \left(\prod_{l=0}^{n-1} w\left(\phi^l\right)\right)B^\phi_{r,0,n}\right\|_{0,K}
	\end{equation*}
	so that (ii) follows.
	
	Conversely, in order to show that (ii) implies (i), we note that, since $\phi$ has stable orbits, for every $j\in\N_0$, the values of $f^{(j)}(\phi^n)$ on $K$ are contained in a compact set which is independent of $n$ and thus can be estimated by a constant. Hence, for arbitrary $f\in\E(X)$, by (\ref{eq:s-th derivative}) and the limits appearing in (ii),
	\begin{equation*}
		\forall\,K\subset X\mbox{ compact}, s\in\N_0:\,\lim_{n\rightarrow\infty}\frac{1}{n}\left\|C_{w,\phi}^n f\right\|_{s,K}=0
	\end{equation*}
	which shows (i).
\end{proof}

\begin{theorem}\label{theo: mean ergodicity on smooth}
	Let $\phi:X\rightarrow X$ and $w:X\rightarrow\C$ be smooth functions with $\{x\in X;\,w(\phi^m(x))\neq 0\}$ being dense in $X$ for every $m\in\N_0$. Consider the following conditions.
	\begin{enumerate}
		\item[(i)] $\phi$ has stable orbits and for every compact set $K\subset X$, $s\in\N_0$, and $h\in\{0,\ldots,s\}$ there holds
		\begin{equation}\label{eq: 1/n condition}
			\lim_{n\rightarrow\infty}\frac{1}{n}\left\|\sum_{r=h}^s\left(\prod_{l=0}^{n-1} w(\phi^l)\right)^{(s-r)}B^\phi_{r,h,n}\right\|_{0,K}=0
		\end{equation}
		as well as
		\begin{equation}\label{eq: sufficient for Cesaro boundedness}
			\sup_{m\in\N}\frac{1}{m}\sum_{n=1}^m \left\|\sum_{r=h}^s \binom{s}{r}\left(\prod_{l=0}^{n-1} w(\phi^l)\right)^{(s-r)}B^\phi_{r,h,n}\right\|_{0,K}<\infty.
		\end{equation}
		\item[(ii)] $C_{w,\phi}$ is (uniformly) mean ergodic on $\E(X)$.
		\item[(iii)] $\phi$ has stable orbits and for every compact set $K\subset X$, $s\in\N_0$, and $h\in\{0,\ldots,s\}$ condition \eqref{eq: 1/n condition} holds.
	\end{enumerate}
	Then (i) implies (ii) and (ii) implies (iii).
\end{theorem}

\begin{proof}
	If (i) holds, it follows from Proposition \ref{prop:Bell to 0} that $\left(\frac{1}{n}C^n_{w,\phi} f\right)_{n\in\N}$ converges to $0$ in $\E(X)$ for every $f\in\E(X)$. Moreover, since $\phi$ has stable orbits, for a fixed compact set $K\subset X$ there is a compact set $L\subset X$ with $\phi^n(K)\subseteq L$ for every $n\in\N$. For arbitrary $f\in\E(X)$ and $s\in\N_0$ it follows from \eqref{eq:s-th derivative} for $x\in K$ and $m\in\N$ that
	\begin{align*}
		\left|\frac{1}{m}\sum_{n=1}^m\left(C^n_{w,\phi} f\right)^{(s)}(x)\right|&\leq \sum_{h=0}^s \frac{1}{m}\sum_{n=1}^m\left|\sum_{r=h}^s\binom{s}{r}\left(\prod_{l=0}^{n-1} w(\phi^l)\right)^{(s-r)}(x)B^\phi_{r,h,n}(x)\right| \times\\
		\times\left|f^{(h)}(\phi^n(x))\right|
		&\leq \left\|f\right\|_{s,L} \sum_{h=0}^s \frac{1}{m}\sum_{n=1}^m\left\|\sum_{r=h}^s\binom{s}{r}\left(\prod_{l=0}^{n-1} w(\phi^l)\right)^{(s-r)}B^\phi_{r,h,n}\right\|_{0,K}
	\end{align*}
	which by \eqref{eq: sufficient for Cesaro boundedness} implies the Ces\`aro boundedness of $C_{w,\phi}$ on $\E(X)$. Since $\E(X)$ is a Montel space, by Theorem \ref{theo:ergodicity and transposed} (b) we conclude that $C_{w,\phi}$ is uniformly mean ergodic.
	
	Next, if $C_{w,\phi}$ is mean ergodic, by Proposition \ref{prop:Bell to 0} (iii) follows.
\end{proof}

\begin{remark}\label{special cases}
	\begin{enumerate}
		\item[(i)] In case $\phi:X\rightarrow X$ is a diffeomorphism it follows from \cite[Proposition 3.9]{K2019p} together with Brouwer's Invariance of Domain Theorem, that $\{x\in X;\,w(\phi^m(x))\neq 0\}$ is dense in $X$ for every $m\in\N_0$ if (and only if) $\{x\in X;\, w(x)\neq 0\}$ is dense in $X$.
		
		\item[(ii)] For the special case of a constant weight $w(x)=\alpha\in\C\backslash\{0\}$, \eqref{eq: 1/n condition} in Theorem \ref{theo: mean ergodicity on smooth} simplifies to
		\begin{enumerate}
			\item[(\ref{eq: 1/n condition}')]
			\begin{equation*}
				\lim_{n\rightarrow\infty}\frac{|\alpha^n|}{n}\left\|B^\phi_{s,h,n}\right\|_{0,K}=0.
			\end{equation*}
		\end{enumerate}
		while \eqref{eq: sufficient for Cesaro boundedness} turns into
		\begin{enumerate}
			\item[(\ref{eq: sufficient for Cesaro boundedness}')]
			\begin{equation*}
				\sup_{m\in\N}\frac{1}{m}\sum_{n=1}^m |\alpha^n| \left\| B^\phi_{s,h,n}\right\|_{0,K}<\infty
			\end{equation*}
		\end{enumerate}
	\end{enumerate}
\end{remark}

\section{Mean ergodic composition operators on $\D'(X)$}\label{sec: mean ergodic on D'}

In this section, for an open subset $X\subseteq\R$, we study mean ergodicity of weighted composition operators on $\D'(X)$ with diffeomorphic symbol, where as usual $\D'(X)$ is equipped with the strong dual topology $\beta(\D'(X),\D(X))$, i.e.\ the topology of uniform convergence on bounded subsets of $\D(X)$. Recall that for a diffeomorphism $\phi:X\rightarrow X$ (or, more generally, a smooth function $\phi:X\rightarrow X$ for which $\phi'$ does not have zeros - so that $\phi$ is injective, in particular) there is a unique continuous linear operator $C_\phi$ on $\D'(X)$ which satisfies $C_\phi f=f\circ\phi$ for every $f\in C(X)$. It holds
\begin{equation*}
	\forall\,u\in \D'(X), \varphi\in\D(X):\,\langle C_\phi u, \varphi\rangle=\left\langle u, \left(\varphi \frac{1}{|\phi'|}\right)\circ\phi^{-1}\right\rangle,
\end{equation*}
where $\langle\cdot,\cdot\rangle$ denotes the duality bracket between $\D'(X)$ and $\D(X)$ (cf.\ \cite[Section 6.1]{HoermanderPDO1}). If additionally $w:X\rightarrow\C$ is smooth we define the weighted composition operator $C_{w,\phi}$ with weight $w$ and symbol $\phi$ as $C_{w,\phi}:=M_w\circ C_\phi$, where $M_w$ denotes the multiplication operator by $w$ on $\D'(X)$, i.e.\ $\langle M_w u,\varphi\rangle=\langle u, w\varphi\rangle$, $u\in\D'(X), \varphi\in\D(X)$. Hence, for $m\in\N_0$ and $u\in\D'(X), \varphi\in\D(X)$
\begin{align*}
	\left\langle C_{w,\phi}^m u,\varphi\right\rangle&=\left\langle u, \left(\varphi\frac{\prod_{l=0}^{m-1}w\left(\phi^l\right)}{|(\phi^m)'|}\right)\circ(\phi^m)^{-1}\right\rangle\\
	&=\left\langle u, \left(\varphi\prod_{l=0}^{m-1}\frac{w\left(\phi^l\right)}{|\phi'\circ\phi^l|}\right)\circ(\phi^m)^{-1}\right\rangle.
\end{align*}
We begin this section with a trivial but important remark.

\begin{remark}\label{rem:trivial but important}
	Let $\phi:X\rightarrow X$ be a diffeomorphism and let $w:X\rightarrow\C$ be smooth. Then
	\begin{equation*}
		\forall\,u\in\D'(X):\,\supp(C_{w,\phi} u)\subseteq\phi\left(\supp u\right).
	\end{equation*}
	In particular, for every distribution with compact support $u$ it follows that $C_{w,\phi}u$ has again compact support, i.e.\ $C_{w,\phi}\left(\E'(X)\right)\subseteq\E'(X)$. Now, we fix $u\in\E'(X)$ and let $\varphi\in\D(X)$ be such that $\varphi=1$ in a neighborhood of $\phi(\supp u)$. Denoting the duality bracket between $\D'(X)$ and $\D(X)$ with $\langle\cdot,\cdot\rangle_{\D',\D}$ and between $\E'(X)$ and $\E(X)$ with $\langle\cdot,\cdot\rangle_{\E',\E}$ for a moment, for $f\in\E(X)$ we have
	\begin{align*}
		\langle C_{w,\phi}u,f\rangle_{\E',\E}&=\langle C_{w,\phi}u,\varphi f\rangle_{\D',\D}=\left\langle u,\left(\left(\frac{w}{|\phi'|} f\right)\circ\phi^{-1}\right)\varphi\circ\phi^{-1}\right \rangle_{\D',\D}\\&=\left\langle u,\left(\frac{w}{|\phi'|} f\right)\circ\phi^{-1}\right \rangle_{\E',\E},
	\end{align*}
	where we have used $\varphi\circ\phi^{-1}=1$ in a neighborhood of $\supp u$ in the last equality. Thus, for the transpose of the restriction of $C_{w,\phi}$ to $\E'(X)$ we have
	\begin{equation*}
		\left(C_{w,\phi|\E'(X)}\right)^t:\E(X)\rightarrow\E(X), f\mapsto \left(\frac{w}{|\phi'|} f\right)\circ\phi^{-1}.
	\end{equation*}
\end{remark}

\begin{proposition}\label{prop:mean ergodic on D' implies mean ergodic on E'}
	Let $\phi:X\rightarrow X$ be a diffeomorphism and $w:X\rightarrow\C$ be smooth such that $\{x\in X;\,w(x)\neq 0\}$ is dense in $X$ and such that $C_{w,\phi}$ is mean ergodic on $\D'(X)$. Then, $C_{w,\phi|\E'(X)}$ is mean ergodic on $\E'(X)$, where the latter is equipped with the strong dual topology $\beta(\E'(X),\E(X))$.
\end{proposition}

Note that $\D'(X)$ as well as $\E'(X)$ are Montel spaces, so by Theorem \ref{theo:ergodicity and transposed} mean ergodicity and uniform mean ergodicity of operators on these spaces are equivalent.

\begin{proof}
	Since $C_{w,\phi}$ is mean ergodic on $\mathscr{D}'(X)$, by Theorem \ref{theo:ergodicity and transposed} (b), it is in particular topologizable. Thus, by \cite[Corollary 2.10]{K2020}, $\phi$ has stable orbits. Therefore, for fixed $u\in\E'(X)$ and $K:=\supp u$ there is a compact $L\subset X$ such that $\phi^n(K)\subseteq L$ for each $n\in\N_0$. In particular, $\supp C_{w,\phi}^n u\subseteq L$ for every $n\in\N$, i.e.\ $C_{w,\phi}^n u \in\E'(L)$ for all $n\in\N$. Since $\E'(L)$ is a closed subspace of $\D'(X)$ and $(C_{w,\phi}^{[n]}u)_{n\in\N}$ converges in $\D'(X)$ by hypothesis, we conclude that $(C_{w,\phi}^{[n]} u)_{n\in\mathbb{N}}$ converges in $\mathscr{E}'(L)$ with respect to the topology $\beta\left(\mathscr{D}'(X),\mathscr{D}(X)\right)$. However, by \cite[Theorem 4.2.1]{Horvath1966}, $\beta(\D'(X),\D(X))$ and $\beta(\E'(X),\E(X))$ induce the same topology on $\E'(L)$ so that $(C_{w,\phi}^{[n]}u)_{n\in\N}$ converges in $\E'(X)$. Since $u\in\E'(X)$ was chosen arbitrarily, the claim follows.
\end{proof}

\begin{proposition}\label{prop:mean ergodicity on D' implies mean ergodicity on E}
	Let $\phi:X\rightarrow X$ be a diffeomorphism and $w:X\rightarrow\C$ be smooth such that $\{x\in X;\,w(x)\neq 0\}$ is dense in $X$ and such that $C_{w,\phi}$ is (uniformly) mean ergodic on $\D'(X)$. Then, the weighted composition operator with weight $w|(\phi^{-1})'|$ and symbol $\phi^{-1}$ on $\E(X)$, $C_{w|(\phi^{-1})'|,\phi^{-1}}$, is (uniformly) mean ergodic.
\end{proposition}

\begin{proof}
	By Proposition \ref{prop:mean ergodic on D' implies mean ergodic on E'}, $C_{w,\phi}$ is (uniformly) mean ergodic on $\E'(X)$. Since $\E(X)$ is a Fr\'echet-Montel space, the (uniform) mean ergodocity of $C_{w|(\phi^{-1})'|,\phi^{-1}}$ on $\E(X)$ follows from Theorem \ref{theo:ergodicity and transposed} and Remark \ref{rem:trivial but important}.
\end{proof}

\begin{theorem}\label{theo: both stable orbits}
	Let $\phi:X\rightarrow X$ be a diffeomorphism and $w:X\rightarrow\C$ be smooth such that $\{x\in X;\,w(x)\neq 0\}$ is dense in $X$ and such that $C_{w,\phi}$ is (uniformly) mean ergodic on $\D'(X)$. Then, $\phi$ and $\phi^{-1}$ have stable orbits.
\end{theorem}

\begin{proof}
	Since $C_{w,\phi}$ is mean ergodic, $C_{w,\phi}$ is topologizable so that $\phi$ has stable orbits by \cite[Corollary 2.10]{K2020}. Moreover, by Proposition \ref{prop:mean ergodicity on D' implies mean ergodicity on E}, $C_{w|(\phi^{-1})'|,\phi^{-1}}$ is mean ergodic on $\E(X)$. Thus, $\phi^{-1}$ has stable orbits by Theorem \ref{theo: mean ergodicity on smooth} together with the fact that with $\{x\in X;\,w(x)\neq 0\}$ being dense in $X$ the same holds for $\{x\in X;\,w(x)|(\phi^{-1})'(x)|\neq 0\}$.
\end{proof}

\begin{example}
	Let $\phi:\R\rightarrow \R$ be the diffeomorphism defined by $\phi(x)=x/2$. Then the operator $C_{\phi}$ is topologizable but it is neither mean ergodic nor power bounded on $\D'(\R)$.
	
	Indeed, clearly $\phi$ has stable orbits so that by \cite[Corollary 2.10]{K2020} $C_{\phi}$ is topologizable. On the other hand, the inverse $\phi^{-1}(x)=2x$ of $\phi$ does not have stable orbits. By Theorem \ref{theo: both stable orbits} the operator $C_{\phi}$ is not mean ergodic. Since $\D'(\R)$ is Montel, in particular semi-reflexive, an application of Corollary \ref{coro: power bounded on semi-reflexive are mean ergodic} shows that $C_{\phi}$ is not power bounded.
\end{example}

The symbol of the operator from the previous example is a real analytic diffeomorphism that allows us to construct an operator $C_{\phi}$ which is topologizable and not mean ergodic on $\D'(\R)$. However, it is not clear whether there is a mean ergodic composition operator on $\mathscr{D}'(\R)$ which is not power bounded. In fact, Theorem \ref{theo: real analytic diffeos} below shows that if such a composition operator exists on $\mathscr{D}'(\R)$ then - recalling that the derivative of the symbol may not have zeros in order to induce a composition operator on $\mathscr{D}'(\R)$ - its symbol cannot be real analytic.

The rest of this section is devoted to prove Theorem \ref{theo: real analytic diffeos}  which characterizes the real analytic diffeomorphisms $\phi$ on an open interval $X\subseteq \R$ for which $C_\phi$ is mean ergodic on $\mathscr{D}'(X)$ and by which this holds precisely when $C_\phi^2=id_{\mathscr{D}'(X)}$.

We denote by $F_\phi$ the set of fixed points of the mapping $\phi:X\rightarrow X$.
\begin{lemma}\label{previous 1}
	Let $\phi:\R\rightarrow \R$ be a diffeomorphism such that $\phi$ and $\phi^{-1}$ have stable orbits. Assume $\phi' >0$, then for each $x\notin F_\phi$ there exist $x_1,x_2\in F_\phi$ such that $x_1<x<x_2$ and
	\begin{equation*}
		]x_1,x_2[\cap F_\phi =\emptyset.
	\end{equation*}
\end{lemma}

\begin{proof}
	Firstly, if there was no fixed point we had $\phi(y)<y$ or $\phi(y)>y$ for all $y\in \R$. In any case, since $\phi'>0$, we obtain that $\phi$ cannot have stable orbits. This is a contradiction and therefore $F_\phi\neq \emptyset$.
	
	Now, we argue by contradiction. Without loss of generality we assume that there is $x_1<x$ such that $]x_1,+\infty[\cap F_\phi =\emptyset$. If $\phi(y)>y$ for all $y>x_1$, then $\phi$ cannot have stable orbits. Finally, if $\phi(y)<y$ for all $y>x_1$, we have that $\phi^{-1}(y)>y$ for all $y>x_1$ thus $\phi^{-1}$ cannot have stable orbits. In any case we obtain a contradiction. 
	
	This completes the proof because $x$ cannot be an accumulation point of $F_\phi$ since the set of fixed points is closed.
\end{proof} 

\begin{remark}\label{special cases moved}
	Let a smooth function $\phi:X\rightarrow X$ be given. For the special weight $w=|\phi'|$, a straightforward calculation gives
	\begin{equation*}
		\forall\,f\in\E(X), n\in\N_0:\,C_{|\phi'|,\phi}^nf=|(\phi^n)'|f(\phi^n),\, \sign\left(\left(\phi^n\right)'\right)=\left(\sign(\phi')\right)^n.
	\end{equation*}
	Denoting a primitive function of $f\in\E(X)$ by $F$, we have for $r, n\in\N_0$
	\begin{align*}
		\left(C_{|\phi'|,\phi}^{n}f\right)^{(r)}&=\Big(\left(f\circ\phi^n\right)\left(\sign\left(\phi'\right)\right)^n (\phi^n)'\Big)^{(r)}\\
		&=\Big(\left(\sign\left(\phi'\right)\right)^n (F\circ \phi^n)'\Big)^{(r)}\\
		&=\left(C_{\sign(\phi'),\phi}^{n}F\right)^{(r+1)}.
	\end{align*}
	Clearly, from the above equality we derive that $C_{|\phi'|,\phi}$ is (uniformly) mean ergodic in $\mathscr{E}(X)$ whenever $C_{\sign (\phi'),\phi}$ is. On the other hand, suppose that $C_{\sign (\phi'),\phi}$ is (uniformly) mean ergodic in $\mathscr{E}(X)$. Then, $\phi$ has stable orbits (see Theorem \ref{theo: mean ergodicity on smooth}) so that the sequence $\big(C_{\sign(\phi'),\phi}^{[n]}F\big)_{n\in\mathbb{N}}$ is bounded with respect to the compact open topology for every $F\in\mathscr{E}(X)$. Additionally, by the fact that $\phi$ has stable orbits, it also holds true that $\big(\frac{1}{n}C^n_{\sign(\phi'),\phi}F\big)_{n\in\mathbb{N}}$ tends to 0 with respect to the compact open topology. Thus, by the above equation $C_{\sign(\phi'),\phi}$ is Cesàro bounded and satisfies that $\lim_{n\rightarrow\infty}\frac{1}{n}C^n_{\sign(\phi'),\phi}=0$  pointwise in $\E(X)$ . Theorem \ref{theo:ergodicity and transposed} (b) yields that $C_{\sign(\phi'),\phi}$ is (uniformly) mean ergodic on $\mathscr{E}(X)$ if the same is true for $C_{|\phi'|,\phi}$.
\end{remark}

\begin{proposition}\label{prop: real analytic diffeos}
	Let $\phi:\R\rightarrow \R$ be a real analytic diffeomorphism. Assume $C_\phi$ is mean ergodic on $\D'(\R)$. Then, $\phi(\phi(x))=\phi^2(x)=x$, $x\in\R$, and consequently $C_\phi^2(u)=u$, $u\in \D'(\R)$.
\end{proposition}
\begin{proof}
	We proceed by contradiction, assuming there is $y\in\R$ such that $\phi^{-2}(y)\neq y$. By Theorem \ref{theo: both stable orbits} for each $K\subset \R$ compact set there is $L\subset\R$ compact such that 
	\[ \bigcup_{n=0}^\infty (\phi^2)^{n}(K)\subseteq \bigcup_{n=0}^\infty \phi^{n}(K)\subseteq L. \]
	Thus we obtain that $\phi^2$ has stable orbits and we can apply the same argument to obtain that $\phi^{-2}$ has stable orbits. Since $\phi$ is a diffeomorphism we have $\phi'(x)\neq 0$ for every $x\in\R$ so that $\phi'>0$ or $\phi'<0$. Thus, by the chain rule $(\phi^2)'>0$, so that by Lemma \ref{previous 1} there exist $x_1,x_2\in F_{\phi^2}$ such that $x_1<y<x_2$ and
	\begin{equation}\label{non-fixed interval}
		]x_1,x_2[\cap F_{\phi^2} =\emptyset.
	\end{equation}
	On the other hand, by Proposition \ref{prop:mean ergodicity on D' implies mean ergodicity on E} and Remark \ref{special cases moved} we obtain that the weighted composition operator $C_{\sign((\phi^{-1})'),\phi^{-1}}:\E(\R)\rightarrow\E(\R)$ is mean ergodic. Since $\E(\R)$ is reflexive we can apply Theorem \ref{mean ergodic}. Observe that Proposition \ref{prop:Bell to 0} implies that 
	\begin{equation*}
		\lim_{n\to \infty} \left| \frac{1}{n} \sign\left((\phi^{-1})'\right)^n \cdot  (\phi^{-n})^{(s)}(x) \right|=\lim_{n\to \infty} \left| \frac{1}{n} (\phi^{-n})^{(s)}(x) \right|=0
	\end{equation*}
	where we have taken $s\geq 1$, $h=1$ and an arbitrary $K=\{x\}\subset \R$ and where we have used that $B_{s,1}(y_1,\ldots,y_s)=y_s$. Denoting $\psi=\phi^{-2}$ we obtain
	\begin{equation}\label{derivatives of 2n iterate}
		\lim_{n\to \infty} \left| \frac{1}{2n} (\psi^{n})^{(s)}(x) \right|=0,
	\end{equation}
	for an arbitrary $x\in\R$.
	
	Auxiliary Claim 1: $\psi'(x_1)=1$ or $\psi'(x_2)=1$.
	
	Indeed, by the chain rule we have that $(\psi^n)'(x)= \left(\psi'(x)\right)^n$ for any fixed point $x$. Then by \eqref{derivatives of 2n iterate} we have 
	\[0<\psi'(x_1)\leq 1 \text{ \quad and \quad } 0<\psi'(x_2)\leq 1.\]
	Suppose that $0<\psi'(x_1)< 1$ and $0<\psi'(x_2)<1$, then there are $y_1,y_2\in ]x_1,x_2[$ such that $\psi(y_1)<y_1$ and $\psi(y_2)>y_2$. Thus by Bolzano's Theorem there is a fixed point between $y_1$ and $y_2$ which contradicts \eqref{non-fixed interval} and proves the Auxiliary Claim 1.
	
	Without loss of generality we can assume $\psi'(x_1)=1$.
	
	By induction one proves $\left(\psi^n\right)''(x_1)=n\cdot\psi''(x_1)$ for all $n\in\N$. Then for $s=2$ and $x_1$ we have by \eqref{derivatives of 2n iterate}
	\[ 0=\lim_{n\to \infty} \left| \frac{1}{2n} (\psi^{n})''(x_1) \right|=\lim_{n\to \infty} \left| \frac{n}{2n} \psi''(x_1)\right|=\frac{1}{2}\left| \psi''(x_1)\right|.\]
	
	Auxiliary Claim 2: Fix $s\geq3$. Assume that $\psi'(x_1)=1$ and $\left(\psi^n\right)^{(j)}(x_1)=0$ hold for all $n\in\N$ and all $2\leq j\leq s-1$. Then $\left(\psi^n\right)^{(s)}(x_1)=n\cdot\psi^{(s)}(x_1)$ for all $n\in\N$.
	
	Indeed, fixing $n\geq2$ and using the original version of Faà di Bruno's Formula for $(\psi\circ \psi^{n-1})^{(s)}$ (see \cite{Johnson02}) we have
	\begin{align*}
		\left(\psi\circ \psi^{n-1}\right)^{(s)}(x_1)=&\sum_{b_1+2b_2+\dots
			+sb_s=s} \frac{s!}{b_1! b_2!\cdots b_s!} \psi^{(b_1+\dots+b_s)}\left(\psi^{n-1}(x_1)\right)\times\\
		&\times \prod_{j=1}^{s}\left( \frac{\left(\psi^{n-1}\right)^{(j)}(x_1)}{j!} \right)^{b_j}.
	\end{align*}
	By the assumptions, given a summand if there is $2\leq j\leq s-1$ with $b_j\neq 0$ then, the summand is automatically $0$. Therefore the non-zero summands must satisfy $b_1+sb_s=s$. By this together with the fact that $x_1$ is a fixed point we obtain
	\begin{align*}
		\left(\psi^n\right)^{(s)}(x_1)&=\psi^{(s)}(x_1) \left(\left(\psi^{n-1}\right)'(x_1)\right)^s + \psi'(x_1)\left(\psi^{n-1}\right)^{(s)}(x_1)\\
		&= \psi^{(s)}(x_1) + \left(\psi^{n-1}\right)^{(s)}(x_1)
	\end{align*}
	because $(\psi^{n-1})'(x_1)=\left(\psi'(x_1)\right)^{n-1}=1$. Auxiliary Claim 2 is now obtained by applying this argument recursively to $\left(\psi^{n-1}\right)^{(s)}(x_1)$.
	
	Under the assumptions of Auxiliary Claim 2 and using \eqref{derivatives of 2n iterate} on $x_1$ we conclude that 
	\[ 0=\lim_{n\to \infty} \left| \frac{n}{2n} \psi^{(s)}(x_1)\right|=\frac{1}{2}\left| \psi^{(s)}(x_1)\right|,\]
	for each $s\geq 2$. To summarize, $x_1$ satisfies $\psi'(x_1)=1$ and $\psi^{(s)}(x_1)=0, s\geq 2$. Because $\psi$ is real analytic as the inverse of a real analytic function, it thus follows $\psi(x)=x$ for every $x\in\R$ which contradicts \eqref{non-fixed interval}. 
\end{proof}

Our next result should be compared to \cite[Theorem 3.8]{FeGaJo18}.

\begin{theorem}\label{theo: real analytic diffeos}
	Let $\phi:X\rightarrow X$ be a real analytic diffeomorphism on the non-empty, open interval $X\subseteq\R$. Then, the following are equivalent.
	\begin{enumerate}
		\item[(i)] $C_\phi:\mathscr{D}'(X)\rightarrow\mathscr{D}'(X)$ is power bounded.
		\item[(ii)] $C_\phi:\mathscr{D}'(X)\rightarrow\mathscr{D}'(X)$ is mean ergodic
		\item[(iii)] $C_\phi:\mathscr{D}'(X)\rightarrow\mathscr{D}'(X)$ is uniformly mean ergodic.
		\item[(iv)] $\phi^2(x)=x$ for each $x\in X$.
		\item[(v)] $C_\phi$ is periodic with period 2.
	\end{enumerate}
\end{theorem}

\begin{proof}
	Let $E$ and $F$ be two lcHs and let $T$ as well as $S$ be continuous linear operators on $E$ and $F$, respectively. Moreover, let $R:E\rightarrow F$ be a continuous linear bijection such that $R^{-1}$ is continuous, too, such that $R\circ T\circ R^{-1}=S$. It is straightforward to show that $T$ is power bounded or (uniformly) mean ergodic if and only if the same applies to $S$.
	
	Moreover, there is a real analytic diffeomorphism $\chi:X\rightarrow\R$ and $C_\chi:\mathscr{D}'(\R)\rightarrow\mathscr{D}'(X), u\mapsto u\circ\chi$ is a continuous linear bijection with $C_\chi^{-1}=C_{\chi^{-1}}$, where for $u\in \mathscr{D}'(\R)$ and $\varphi\in\mathscr{D}(X)$ as usual $\langle C_\chi u,\varphi\rangle=\left\langle u,\left(\frac{\varphi}{|\chi'|}\right)\circ\chi^{-1}\right\rangle$. Clearly, $C_\chi^{-1}\circ C_\phi\circ C_\chi=C_{\chi\circ\phi\circ\chi^{-1}}$ and the real analytic diffeomorphism $\chi\circ\phi\circ\chi^{-1}$ on $\R$ satisfies $\left(\chi\circ\phi\circ\chi^{-1}\right)^2(x)=x$ for every $x\in\R$ if and only if $\phi^2(x)=x$ for each $x\in X$. Therefore, without loss of generality we can assume $X=\R$.
	
	Next, since $\mathscr{D}'(\R)$ is a Montel space it is in particular a semi-reflexive space. Then by Corollary \ref{coro: power bounded on semi-reflexive are mean ergodic} every power bounded operator $T$ on $\mathscr{D}'(\R)$ is mean ergodic. An application of Theorem \ref{theo:ergodicity and transposed} (a) gives that $T$ is also uniformly mean ergodic. Hence, (ii) follows from (i), and (ii) and (iii) are equivalent. By Proposition \ref{prop: real analytic diffeos}, (ii) implies (iv). Finally, (iv) trivially implies (v) which in turn immediately implies (i).
\end{proof}

\begin{corollary}
	Let $\phi:X\rightarrow X$ be a real analytic diffeomorphism on an open interval $X\subseteq\R$ such that there is $p\in \N$ with $\phi^p(x)=x$ for all $x\in X$. Then $\phi^2(x)=x$ for each $x\in X$.
\end{corollary}

\begin{proof}
	By hypothesis $C_\phi$ is a periodic operator on $\mathscr{D}'(X)$ so it is in particular power bounded. Thus, the claim follows from Theorem \ref{theo: real analytic diffeos}.
\end{proof}

We give examples of real analytic diffeomorphisms $\phi$ satisfying the assumptions of the previous result. In case $\phi'>0$ the unique possibility is the identity map. However, the case $\phi'<0$ is richer. The following example should be compared to \cite[Proposition 3.6 and Example 1]{FeGaJo18}. 

\begin{example}
	Let $\phi:X\rightarrow X$ be a real analytic diffeomorphism on an open interval $X\subseteq\R$ such that $\phi^2(x)=x$ for each $x\in X$.
	
	Assuming that our real analytic diffeomorphism $\phi:X\rightarrow X$ is increasing, it follows that $\phi(x)=x$, for all $x\in X$. Indeed, if there is $y\neq \phi(y)$ we may assume without loss of generality that $y<\phi(y)$ (for otherwise, we consider $\phi^{-1}$). Now, since $\phi$ is increasing we obtain that $y<\phi(y)<\phi^2(y)$ which contradicts the assumption.
	
	Contrary to the case of increasing real analytic diffeomorphisms, if we assume $\phi'<0$, apart from the obvious example $\phi(x)=-(x+c)$ where $c\in\R$, we obtain a whole class of examples as in \cite[Example 1]{FeGaJo18}. For this, let $f:\R\rightarrow\R$ be an even real analytic function such that $|f'(x)|\leq a < 1$ for all $x\in \R$. Then, except for refering to the Implicit Function Theorem in the real analytic class, see e.g.\ \cite[Theorem 2.3.5]{KrantzParks}, a verbatim repetition of the arguments presented in \cite[Example 1]{FeGaJo18}, the equation 
	\[ x+y=f(x-y)\]
	defines a decreasing real analytic symbol $y=\phi(x)$, $x\in\R$, with $\phi^2(x)=x$. Thus, $\phi$ is bijective, hence a real analytic diffeomorphism.
	
	For a concrete example, we consider the function $f(x)=\sqrt{\frac{x^2}{2}+1}$ which yields the real analytic decreasing diffeomorphism 
	\[ \phi:\R\rightarrow\R, x\mapsto -3x + \sqrt{8x^2+2}. \]
	It is easy to see that $\phi=\phi^{-1}$.	
\end{example}

\noindent \textbf{Acknowledgements.} The authors are very grateful to J.\ Bonet from Universitat Polit\`ecnica de Val\`encia for many stimulating discussions on the topic of this paper and cordially acknowledge the permission to include Corollary 2.2 and Theorem 2.3 from a personal communication, and for pointing out \cite[Example 1]{FeGaJo18}. Moreover, the authors thank the anonymous referee for valuable suggestions which helped to improve the presentation of the paper. The research of the second author was partially supported by the project GV Prometeo 2017/102.

\bibliographystyle{abbrv}
\bibliography{bib_KS}

\begin{thebibliography}{10}

\bibitem{ABR2009}
A.~A. Albanese, J.~Bonet, and W.~J. Ricker.
\newblock Mean ergodic operators in {F}r\'{e}chet spaces.
\newblock {\em Ann. Acad. Sci. Fenn. Math.}, 34(2):401--436, 2009.

\bibitem{ABR2010}
A.~A. Albanese, J.~Bonet, and W.~J. Ricker.
\newblock On mean ergodic operators.
\newblock In {\em Vector measures, integration and related topics}, volume 201
  of {\em Oper. Theory Adv. Appl.}, pages 1--20. Birkh\"{a}user Verlag, Basel,
  2010.

\bibitem{Beltran20}
M.~J. Beltr\'{a}n-Meneu.
\newblock Dynamics of weighted composition operators on weighted {B}anach
  spaces of entire functions.
\newblock {\em J. Math. Anal. Appl.}, 492(1):124422, 16, 2020.

\bibitem{BeGCJoJo16}
M.~J. Beltr\'{a}n-Meneu, M.~C. G\'{o}mez-Collado, E.~Jord\'{a}, and D.~Jornet.
\newblock Mean ergodic composition operators on {B}anach spaces of holomorphic
  functions.
\newblock {\em J. Funct. Anal.}, 270(12):4369--4385, 2016.

\bibitem{BGJJ2016mw}
M.~J. Beltr\'{a}n-Meneu, M.~C. G\'{o}mez-Collado, E.~Jord\'{a}, and D.~Jornet.
\newblock Mean ergodicity of weighted composition operators on spaces of
  holomorphic functions.
\newblock {\em J. Math. Anal. Appl.}, 444(2):1640--1651, 2016.

\bibitem{BoPaRi11}
J.~Bonet, B.~de~Pagter, and W.~J. Ricker.
\newblock Mean ergodic operators and reflexive {F}r\'{e}chet lattices.
\newblock {\em Proc. Roy. Soc. Edinburgh Sect. A}, 141(5):897--920, 2011.

\bibitem{BoDo2011A}
J.~Bonet and P.~Doma\'{n}ski.
\newblock A note on mean ergodic composition operators on spaces of holomorphic
  functions.
\newblock {\em Rev. R. Acad. Cienc. Exactas F\'{\i}s. Nat. Ser. A Mat. RACSAM},
  105(2):389--396, 2011.

\bibitem{BoDo11B}
J.~Bonet and P.~Doma\'{n}ski.
\newblock Power bounded composition operators on spaces of analytic functions.
\newblock {\em Collect. Math.}, 62(1):69--83, 2011.

\bibitem{FeGaJo18}
C.~Fern\'{a}ndez, A.~Galbis, and E.~Jord\'{a}.
\newblock Dynamics and spectra of composition operators on the {S}chwartz
  space.
\newblock {\em J. Funct. Anal.}, 274(12):3503--3530, 2018.

\bibitem{GCJoJo16}
M.~C. G\'{o}mez-Collado, E.~Jord\'{a}, and D.~Jornet.
\newblock Power bounded composition operators on spaces of meromorphic
  functions.
\newblock {\em Topology Appl.}, 203:141--146, 2016.

\bibitem{HaShZh19}
S.-A. Han and Z.-H. Zhou.
\newblock Mean ergodicity of composition operators on {H}ardy space.
\newblock {\em Proc. Indian Acad. Sci. Math. Sci.}, 129(4):Paper No. 45, 10,
  2019.

\bibitem{HoermanderPDO1}
L.~H{\"o}rmander.
\newblock {\em The analysis of linear partial differential operators. {I}}.
\newblock Classics in Mathematics. Springer-Verlag, Berlin, 2003.
\newblock Distribution theory and Fourier analysis, Reprint of the second
  (1990) edition [Springer, Berlin; MR1065993 (91m:35001a)].

\bibitem{Horvath1966}
J.~Horv\'{a}th.
\newblock {\em Topological vector spaces and distributions. {V}ol. {I}}.
\newblock Addison-Wesley Publishing Co., Reading, Mass.-London-Don Mills, Ont.,
  1966.

\bibitem{Johnson02}
W.~P. Johnson.
\newblock The curious history of {F}a\`a di {B}runo's formula.
\newblock {\em Amer. Math. Monthly}, 109(3):217--234, 2002.

\bibitem{JoRA20}
E.~Jord\'{a} and A.~Rodr\'{\i}guez-Arenas.
\newblock Ergodic properties of composition operators on {B}anach spaces of
  analytic functions.
\newblock {\em J. Math. Anal. Appl.}, 486(1):123891, 14, 2020.

\bibitem{JoSaSP20}
D.~Jornet, D.~Santacreu, and P.~Sevilla-Peris.
\newblock Mean ergodic composition operators in spaces of homogeneous
  polynomials.
\newblock {\em J. Math. Anal. Appl.}, 483(1):123582, 11, 2020.

\bibitem{JoSaSP21}
D.~Jornet, D.~Santacreu, and P.~Sevilla-Peris.
\newblock Mean ergodic composition operators on spaces of holomorphic functions
  on a {B}anach space.
\newblock {\em J. Math. Anal. Appl.}, 500(2):125139, 16, 2021.

\bibitem{K2019p}
T.~Kalmes.
\newblock Power bounded weighted composition operators on function spaces
  defined by local properties.
\newblock {\em J. Math. Anal. Appl.}, 471(1-2):211--238, 2019.

\bibitem{K2020}
T.~Kalmes.
\newblock Topologizable and power bounded weighted composition operators on
  spaces of distributions.
\newblock {\em Ann. Polon. Math.}, 125(2):139--154, 2020.

\bibitem{KrantzParks}
S.~G. Krantz and H.~R. Parks.
\newblock {\em A primer of real analytic functions}, volume~4 of {\em Basler
  Lehrb\"{u}cher [Basel Textbooks]}.
\newblock Birkh\"{a}user Verlag, Basel, 1992.

\bibitem{Krengel}
U.~Krengel.
\newblock {\em Ergodic theorems}, volume~6 of {\em De Gruyter Studies in
  Mathematics}.
\newblock Walter de Gruyter \& Co., Berlin, 1985.
\newblock With a supplement by Antoine Brunel.

\bibitem{SeMeBo20}
W.~Seyoum, T.~Mengestie, and J.~Bonet.
\newblock Mean ergodic composition operators on generalized {F}ock spaces.
\newblock {\em Rev. R. Acad. Cienc. Exactas F\'{\i}s. Nat. Ser. A Mat. RACSAM},
  114(1):Paper No. 6, 11, 2020.

\end{thebibliography}

\end{document}